\RequirePackage{ifpdf}
\ifpdf 
\documentclass[pdftex]{sigma}
\else
\documentclass{sigma}
\fi

\usepackage{dsfont}

\newcommand{\reals}{\mathbb{R}}

\newcommand{\complex}{\mathbb{C}}

\newcommand{\bracketa}[1]{\big[#1\big]}
\newcommand{\bracketb}[1]{\Big[#1\Big]}

\newcommand{\angles}[1]{\left\langle #1 \right\rangle}
\newcommand{\pb}[1]{\left\{#1\right\}}
\newcommand{\com}[1]{\left[#1\right]}

\newcommand{\Ccom}[3]{\big[[#1,#2],#3\big]}

\newcommand{\paraa}[1]{\big(#1\big)}

\newcommand{\diag}{\operatorname{diag}}

\newcommand{\spacearound}[1]{\quad#1\quad}
\newcommand{\equivalent}{\spacearound{\Longleftrightarrow}}

\renewcommand{\mid}{\mathds{1}}
\newcommand{\X}{\mathcal{X}}
\newcommand{\Y}{\mathcal{Y}}
\newcommand{\A}{\mathcal{A}}
\newcommand{\B}{\mathcal{B}}
\newcommand{\g}{\mathfrak{g}}
\newcommand{\h}{\mathfrak{h}}

\newcommand{\sumjd}{\sum_{j=1}^d}
\newcommand{\DMSA}{\textsf{DMSA}}

\renewcommand{\sl}{\mathfrak{sl}}
\newcommand{\su}{\mathfrak{su}}
\newcommand{\DX}{\Delta_\X}
\newcommand{\U}{\mathfrak{U}}
\newcommand{\Th}{T^{\hbar}}
\newcommand{\Td}{T^\dagger}

\newcommand{\xv}{\vec{x}}
\renewcommand{\L}{\Lambda}
\newcommand{\Ld}{\L^\dagger}

\newcommand{\vspan}[1]{\langle#1\rangle}
\newcommand{\ep}{e^{+}}
\newcommand{\emi}{e^{-}}
\newcommand{\ea}{e_{\alpha}}
\newcommand{\ema}{e_{-\alpha}}

\newcommand{\eap}{\ep_{\alpha}}
\newcommand{\eam}{\emi_{\alpha}}
\newcommand{\ebp}{\ep_{\beta}}
\newcommand{\ebm}{\emi_{\beta}}
\newcommand{\ha}{h_{\alpha}}
\newcommand{\hb}{h_{\beta}}

\newcommand{\ad}{\operatorname{ad}}
\newcommand{\Spec}{\operatorname{Spec}}
\newcommand{\phitt}{\phi_{\theta\theta'}}

\newcommand{\alphat}{\tilde{\alpha}}

\newcommand{\Dt}{\tilde{D}}
\newcommand{\Dtd}{\tilde{D}^\dagger}
\newcommand{\dt}{\tilde{d}}

\newcommand{\twovec}[2]{\begin{pmatrix}#1 \\ #2\end{pmatrix}}
\newcommand{\twomatrix}[4]{\begin{pmatrix} #1 & #2 \\ #3 & #4\end{pmatrix} }
\newcommand{\vphi}{\varphi}
\newcommand{\vphit}{\tilde{\varphi}}
\newcommand{\Cl}{\mathit{Cl}}

\begin{document}

\allowdisplaybreaks

\renewcommand{\thefootnote}{$\star$}

\renewcommand{\PaperNumber}{042}

\FirstPageHeading

\ShortArticleName{Discrete Minimal Surface Algebras}

\ArticleName{Discrete Minimal Surface Algebras\footnote{This paper is a
contribution to the Special Issue ``Noncommutative Spaces and Fields''. The
full collection is available at
\href{http://www.emis.de/journals/SIGMA/noncommutative.html}{http://www.emis.de/journals/SIGMA/noncommutative.html}}}

\Author{Joakim ARNLIND~$^\dag$ and Jens HOPPE~$^\ddag$}

\AuthorNameForHeading{J.~Arnlind and J.~Hoppe}

\Address{$^\dag$~Institut des Hautes \'Etudes Scientif\/iques, Le Bois-Marie, 35, Route de Chartres,\\
 \hphantom{$^\dag$}~91440 Bures-sur-Yvette, France}
\EmailD{\href{mailto:arnlind@ihes.fr}{arnlind@ihes.fr}}

\Address{$^\ddag$~Eidgen\"ossische Technische Hochschule, 8093 Z\"urich, Switzerland\\
\hphantom{$^\ddag$}~(on leave of absence from  Kungliga Tekniska H\"ogskolan, 100 44 Stockholm, Sweden)}
\EmailD{\href{mailto:hoppe@itp.phys.ethz.ch}{hoppe@itp.phys.ethz.ch}, \href{mailto:hoppe@math.kth.se}{hoppe@math.kth.se}}

\ArticleDates{Received March 23, 2010, in f\/inal form May 18, 2010;  Published online May 26, 2010}

\Abstract{We consider discrete minimal surface algebras (\DMSA) as
  generalized non\-com\-mu\-ta\-ti\-ve analogues of minimal surfaces in
  higher dimensional spheres. These algebras appear naturally in
  membrane theory, where sequences of their representations are used
  as a regularization. After showing that the def\/ining relations of
  the algebra are consistent, and that one can compute a basis of the
  enveloping algebra, we give several explicit examples of
  \DMSA s in terms of subsets of $\sl_n$ (any semi-simple Lie algebra
  providing a trivial example by itself). A special class of \DMSA s
  are Yang--Mills algebras.
  The representation graph is introduced to study representations of
  \DMSA s of dimension $d\leq 4$, and properties of representations
  are related to properties of graphs. The representation graph of a
  tensor product is (generically) the Cartesian product of the
  corresponding graphs. We provide explicit examples of irreducible
  representations and, for coinciding eigenvalues, classify all the
  unitary representations of the corresponding algebras.}


\Classification{81R10; 06B15}

\section{Introduction}

 Noncommutative analogues of manifolds have been studied in
many dif\/ferent contexts. One way of constructing such objects is to
relate the Poisson structure of a manifold $M$ to the commutator
structure of a sequence of matrix algebras $\A_\hbar$ (parametrized by
$\hbar>0$), where the dimension of the matrices increases as $\hbar\to
0$. Namely, for some set of values of the parameter $\hbar$, one
def\/ines a map $\Th:C^{\infty}(M)\to\A_\hbar$ such that
\begin{gather}\label{eq:ToeplizCondition}
  \lim_{\hbar\to 0}\left|\left|\Th\paraa{\{f,g\}}-\frac{1}{i\hbar}\bracketa{\Th(f),\Th(g)}\right|\right|=0
\end{gather}
for all $f,g\in C^{\infty}(M)$.

For surfaces, the map $\Th$ has been constructed in several dif\/ferent
ways. One rather concrete approach is to consider the surface
$\Sigma$ as embedded in an ambient manifold $M$ with embedding
coordinates $x_1(\sigma_1,\sigma_2),\ldots,x_d(\sigma_1,\sigma_2)$. If
the Poisson brackets satisfy
\begin{gather*}
  \pb{x_i,x_j} = p_{ij}(x_1,\ldots,x_d)
\end{gather*}
where $p_{ij}(x_1,\ldots,x_d)$ are polynomials, then one def\/ines an
algebra of non-commuting variables $X_1,\ldots,X_d$ such that the
following relations hold
\begin{gather*}
  \bracketa{X_i,X_j} = i\hbar \Psi(p_{ij})(X_1,\ldots,X_d)
\end{gather*}
where $\Psi$ is an ordering map, mapping commutative polynomials to
non-commutative ones (such that when composed with the
projection back to commutative polynomials one gets the identity
map). Thus, for any pair of polynomials $p$, $q$ it holds that
\begin{gather*}
  \Psi(\{p,q\})-\frac{1}{i\hbar}[\Psi(p),\Psi(q)] = O(\hbar),
\end{gather*}
and by considering representations of the above def\/ined algebra one
obtains a sequence of matrix algebras and maps $\Th$ such that relation
\eqref{eq:ToeplizCondition} holds for all polynomial functions (see
\cite{a:phdthesis} for details). It is natural to demand that a notion
of noncommutative analogues of manifolds should have some features
that can be traced back to the geometry of the original manifold. For
surfaces, the genus is the obvious invariant, and one can show that
the above procedure gives rise to algebras whose representation theory
encodes geometric data \cite{abhhs:noncommutative}.

For a surface $\Sigma$ embedded in $\reals^d$, the Laplace--Beltrami
operator on $\Sigma$ acting on the embedding coordinates
$x_1(\sigma_1,\sigma_2),\ldots,x_d(\sigma_1,\sigma_2)$ can be written
as
\begin{gather*}
  \Delta(x_i) = \sumjd\pb{\pb{x_i,x_j},x_j}
\end{gather*}
where
$\pb{f,h}(\sigma_1,\sigma_2)=\frac{1}{\sqrt{g}}(\partial_1f\partial_2h-\partial_1h\partial_2f)$
and $g$ is the determinant of the induced metric on $\Sigma$. With
this notation, minimal surfaces in $S^{d-1}$ can be found by
constructing embedding coordinates such that
\begin{gather*}
  \sumjd\pb{\pb{x_i,x_j},x_j} = -2x_i
\end{gather*}
subject to the constraint $\sumjd x_j^2=1$.  In the above spirit of
replacing Poisson brackets by commutators, corresponding
noncommutative minimal surfaces are def\/ined by the relations
\begin{align}
  &\hbar^2\sumjd\Ccom{X_i}{X_j}{X_j} = 2X_i.\label{eq:doubleCommutators}
\end{align}

 Another example where equations like \eqref{eq:doubleCommutators}
arise is in the context of a physical theory of ``Membranes'' \cite{h:phdthesis}. The
equations of motion for a membrane moving in $d+1$ dimensional
Minkowski space (with a particular choice of coordinates) can be
written as
\begin{gather*}
  \partial^2_tx_i = \sumjd\pb{\pb{x_i,x_j},x_j}\qquad\text{ and }\qquad
  \sumjd\pb{\partial_tx_j,x_j} = 0.
\end{gather*}
A regularized theory is given by $d$ time-dependent Hermitian
matrices $X_i$ satisfying the equations
\begin{gather}\label{eq:membraneEqMotion}
  \partial^2_tX_i=-\hbar^2\sumjd\Ccom{X_i}{X_j}{X_j}
  \quad\text{ and }\quad\sumjd\com{\partial_tX_i,X_j}=0.
\end{gather}
In the f\/irst equation, we can separate time from the matrices by
making the ansatz
\begin{gather}\label{eq:timeSeparationAnsatz}
  X_i(t) = \sum_{j=1}^d\frac{a}{\hbar}\paraa{e^{A(at+b)}}_{ij}M_j,
\end{gather}
(see also \cite{h:someClassicalSolutions,aht:spinning}) where $A$ is a
$d\times d$ antisymmetric matrix such that
$A^2=\diag(-\mu_1,\ldots,-\mu_d)$ with
$\mu_{2i-1}=\mu_{2i}$ for $i=1,2,\ldots,\lfloor
d/2\rfloor$. Then $X_i$ (as def\/ined in~\eqref{eq:timeSeparationAnsatz}) solves the f\/irst equation in~\eqref{eq:membraneEqMotion} provided
\begin{gather}\label{eq:MDoubleCommutator}
  \sumjd\Ccom{M_i}{M_j}{M_j} = \mu_iM_i.
\end{gather}
Motivated by the above examples, we set out to study algebras
generated by relations \eqref{eq:MDoubleCommutator} (for \emph{arbitrary}
$\mu_i$'s). In particular, we shall study their representation
theory for $d\leq 4$.

In Section \ref{sec:dmsa} we introduce Discrete Minimal Surface
Algebras, and after showing that a~basis of the enveloping
algebra can be computed via the Diamond lemma, some properties are
investigated and several examples are given. Section
\ref{sec:hermitianReps} deals with Hermitian representations of the
f\/irst non-trivial algebras, for which the representation graph is
introduced, and properties of representations are related to
properties of graphs. Finally, we provide explicit examples of
irreducible representations and their corresponding graphs and, in the
case when $\Spec(\A)=\{\mu\}$, we classify all unitary representations.

\section{Discrete minimal surface algebras}\label{sec:dmsa}

 Let $\g$ be a Lie algebra over $\complex$. For any f\/inite
subset $\X=\{x_1,x_2,\ldots,x_d\}\subseteq\g$ we def\/ine a linear map
$\DX:\g\to\g$ by
\begin{align*}
  \DX(a) = \sumjd\Ccom{a}{x_j}{x_j},
\end{align*}
and let $\vspan{\X}$ denote the vector space spanned by
the elements in $\X$.

\begin{definition}[\DMSA]
  Let $\g$ be a Lie algebra and let $\X=\{x_1,x_2,\ldots,x_d\}$ be a
  set of linearly independent elements of $\g$. We call
  $\A=(\g,\X)$ a \emph{discrete minimal surface algebra} (\DMSA) (of
  dimension $d$) if there exist complex numbers
  $\mu_1,\ldots,\mu_d$ such that $\DX(x_i)=\mu_ix_i$ for
  $i=1,\ldots,d$. The set $\{\mu_1,\mu_2,\ldots,\mu_d\}$
  is called \emph{the spectrum of $\A$} and is denoted by $\Spec(\A)$.
\end{definition}

\begin{definition}
  Two \DMSA s $\A=(\g,\X)$ and $\B=(\g',\Y)$ are \emph{isomorphic} if there
  exists a Lie algebra homomorphism $\phi:\g\to\g'$ such that
  $\phi|_{\vspan{\X}}$ is a vector space isomorphism of $\vspan{\X}$
  onto~$\vspan{\Y}$.
\end{definition}

 Note that $\Spec(\A)$ is not an invariant within an
isomorphism class. Let $\g$ be a Lie algebra with structure constants
$f_{ij}^k$ (relative to the basis $x_1,\ldots,x_n$) and let $\g'$ be a
Lie algebra with structure constants $cf_{ij}^k$ (relative to the
basis $y_1,\ldots,y_n$), for some non-zero complex number $c$.  Then
the two \DMSA s $\A=(\g,\X=\{x_1,\ldots,x_d\})$ and
$\B=(\g',\Y=\{y_1,\ldots,y_d\})$ will be isomorphic (through
$\phi(x_i)=y_i/c$), and if $\Spec(\A)=\{\mu_1,\ldots,\mu_d\}$
then $\Spec(\B)=\{c^2\mu_1,\ldots,c^2\mu_d\}$.

\begin{definition}
  Let $\A=(\g,\X)$ be a \DMSA\ with
  $\Spec(\A)=\{\mu_1,\ldots,\mu_d\}$ and let
  $\complex\angles{\X}$ denote the free associative algebra over
  $\complex$, generated by the set
  $\mathcal{X}=\{x_1,x_2,\ldots,x_d\}$.  The enveloping
  algebra of $\A$ is def\/ined as the quotient
  $\U_d(\A)=\complex\angles{\X}\slash \angles{\DX(x_i)-\mu_ix_i}$,
  where $[a,b]=ab-ba$.
\end{definition}

\begin{remark}
  Relating the above enveloping algebra to minimal surfaces in
    $S^{d-1}$, one could impose the relation
    $x_1^2+x_2^2+\cdots+x_d^2=\nu\mid$ together with $\DX(x_i)=\mu
    x_i$. Combining these relations gives the equations
  \begin{align*}
    \sumjd x_jx_ix_j = \frac{1}{2}(2\nu-\mu)x_i
  \end{align*}
  for $i=1,\ldots,d$.
\end{remark}

 The following result shows how to compute a basis of the
enveloping algebra. Note that in~\cite{bdv:inhomyangmills}, it was
proven that a class of so called inhomogeneous Yang--Mills algebras
(including $\U_d(\A)$) has the Poincar\'e--Birkhof\/f--Witt property.

\begin{proposition}
  A basis of $\U_d(\A)$ is provided by the set of words on
  $\{x_1,\ldots,x_d\}$ that do not contain any of the following
  sub-words
  \begin{align*}
    x_d^2x_1,\,x_d^2x_2,\,\ldots,\,x_d^2x_{d-1},\,x_dx_{d-1}^2.
  \end{align*}
\end{proposition}

\begin{proof}
  In the following proof we will use the notation and terminology of
  \cite{b:diamondLemma}, to which we refer for details. The def\/ining relations of the algebra
  \begin{gather*}
    \sumjd\Ccom{x_i}{x_j}{x_j} = \mu_ix_i
  \end{gather*}
  for $i=1,\ldots,d$ will be put into the following reduction system $S$
  \begin{gather*}
     \sigma_i = (W_i,f_i) = \bigg(x_d^2x_i,2x_dx_ix_d-x_ix_d^2+\lambda_ix_i-\sum_{j=1}^{d-1}\Ccom{x_i}{x_j}{x_j}\bigg),
    \qquad 1\leq i<d,\\
     \sigma_d = (W_d,f_d) = \bigg(x_dx_{d-1}^2,2x_{d-1}x_dx_{d-1}-x_{d-1}^2x_d+\lambda_dx_d
      -\sum_{j=1}^{d-2}\Ccom{x_d}{x_j}{x_j}\bigg).
  \end{gather*}
  Next, we need to def\/ine a semi-group partial ordering on words, that
  is compatible with $S$, i.e.\ every word in $f_i$ should be less than
  $W_i$.  Let $w_1$, $w_2$ be two words on $x_1,\ldots,x_d$; if $w_1$ has
  smaller length than $w_2$ then we set $w_1<w_2$. If $w_1$ and $w_2$
  have the same length, we set $w_1<w_2$ if $w_1$ precedes $w_2$
  lexicographically, where the lexicographical ordering is induced by
  $x_1<x_2<\cdots<x_d$. These def\/initions imply that if $w_1<w_2$ then
  $aw_1b<aw_2b$ for any words $a$, $b$ (which def\/ines a semi-group
  partial ordering). It is also easy to check that this ordering is
  compatible with $S$. Does the ordering satisfy the descending chain
  condition? Let $w_1\geq w_2\geq w_3 \geq\cdots$ be an inf\/inite
  sequence of decreasing words. Clearly, since the length of $w_i$ is
  a positive integer, it must eventually become constant. Thus, for
  all $i>N$ the length of $w_i$ is the same. But this implies that the
  series eventually become constant because there is only a f\/inite
  number of words preceding a given word lexicographically. Hence,
  the ordering satisf\/ies the descending chain condition.

  Now, we are ready to apply the Diamond lemma. If we can show that
  all ambiguities in $S$ are resolvable, then a basis for the algebra
  is provided by the irreducible words. In this case there is only one
  ambiguity in the reduction system $S$. Namely, there are two ways of
  reducing the word $x_d^2x_{d-1}^2$: Either we write it as
  $x_d(x_dx_{d-1}^2)$ and apply $\sigma_d$ or we write it as
  $(x_d^2x_{d-1})x_{d-1}$ and apply $\sigma_{d-1}$. Let us prove that
  $A\equiv x_df_d-f_{d-1}x_{d-1}=0$, i.e.\ the ambiguity is resolvable
  \begin{gather*}
    A =  -x_{d-1}\overbrace{x_d^2x_{d-1}}^{\to f_{d-1}}+\overbrace{x_dx_{d-1}^2}^{\to f_d}x_d
    + \lambda_{d-1}x_{d-1}^2-\lambda_dx_d^2\\
    \phantom{A=}{} -\sum_{j=1}^{d-1}\Ccom{x_{d-1}}{x_j}{x_j}x_{d-1}
    +x_d\sum_{j=1}^{d-2}\Ccom{x_d}{x_j}{x_j}\\
   \phantom{A} =   \sum_{j=1}^{d-2}\Big[x_{d-1},\Ccom{x_{d-1}}{x_j}{x_j}\Big]
    +\sum_{j=1}^{d-2}\Big[x_{d},\Ccom{x_{d}}{x_j}{x_j}\Big].
  \end{gather*}
  Now, let us rewrite the second commutator above as follows:
  \begin{gather*}
    \Big[x_{d},\Ccom{x_{d}}{x_j}{x_j}\Big] =
     (\overbrace{x_d^2x_j}^{\to f_j})x_j-2x_dx_jx_dx_j-x_j^2x_d^2+2x_jx_dx_jx_d\\
  \hphantom{\Big[x_{d},\Ccom{x_{d}}{x_j}{x_j}\Big] }{}  =  2x_jx_dx_jx_d-x_j^2x_d^2+\lambda_jx_j^2-x_j\overbrace{x_d^2x_j}^{\to f_j}
    -\sum_{k=1}^{d-1}\Ccom{x_j}{x_k}{x_k}x_j\\
   \hphantom{\Big[x_{d},\Ccom{x_{d}}{x_j}{x_j}\Big] }{} =   \sum_{k=1}^{d-1}\Big[x_{j},\Ccom{x_{j}}{x_k}{x_k}\Big].
  \end{gather*}
  By introducing the sum, we obtain
  \begin{gather*}
    \sum_{j=1}^{d-2}\Big[x_{d},\Ccom{x_{d}}{x_j}{x_j}\Big] =
     -\sum_{j=1}^{d-2}\Big[x_{d-1},\Ccom{x_{d-1}}{x_j}{x_{j}}\Big]
    +\sum_{j,k=1}^{d-2}\Big[x_{j},\Ccom{x_{j}}{x_k}{x_k}\Big]\\
  \hphantom{\sum_{j=1}^{d-2}\Big[x_{d},\Ccom{x_{d}}{x_j}{x_j}\Big]}{}
    =  -\sum_{j=1}^{d-2}\Big[x_{d-1},\Ccom{x_{d-1}}{x_j}{x_{j}}\Big],
  \end{gather*}
  which implies that $A=0$. From the Diamond lemma, we can now
  conclude that the set of all irreducible words, with respect to the
  reduction system $S$, provides a basis of the algebra.
\end{proof}

  Let us start by noting that any semi-simple Lie algebra $\g$ is itself
a \DMSA. Namely, if we let $\X=\{x_1,x_2,\ldots,x_d\}$ be a basis of
$\g$ such that $K(x_i,x_j)=\delta_{ij}$ (where $K$ denotes the Killing
form), then the structure constants will be totally anti-symmetric
which implies that
\begin{gather*}
  \DX(x_i)=\sumjd\Ccom{x_i}{x_j}{x_j} = \sum_{j,k,l=1}^df_{ij}^kf_{kj}^lx_l=
  \sum_{l=1}^dK(x_i,x_l)x_l = x_i.
\end{gather*}
Thus, in such a basis $(\g,\g)$ is a \DMSA\ for any semi-simple Lie
algebra $\g$. On the other hand, if $\g$ is nilpotent, it follows that
the map $\DX$ is nilpotent. Thus, for a \DMSA\ related to a~nilpotent
Lie algebra, it must hold that $\Spec(\A)=\{0\}$. Note that the class
of \DMSA\ for which $\Spec{\A}=\{0\}$ has also been studied under the
name of \emph{$($Lie$)$ Yang--Mills algebras}
\cite{n:lecturenoncgauge,cdv:yangMillsAlgebra}, and their
representation theory has been studied in
\cite{hs:repYangMillsAlgebras}.

Apart from semi-simple Lie algebras, it is also true that Clif\/ford
algebras satisfy the relations imposed in the enveloping
algebra. Let $\Cl_{p,q}$ be a Clif\/ford algebra generated by
$e_1,\ldots,e_d$ (with $d=p+q$), satisfying
\begin{alignat*}{3}
  & e_i^2=1\qquad&&\text{for} \ \ i=1,\ldots,p, & \\
  & e_i^2=-1\qquad&&\text{for} \ \ i=p+1,\ldots,p+q, & \\
  & e_ie_j=-e_je_i\qquad&& \text{when} \ \ i\neq j. &
\end{alignat*}
It is then easy to check that
\begin{align*}
  \sumjd\Ccom{e_i}{e_j}{e_j} =
  \begin{cases}
    4(p-q-1)e_i &\text{if} \ \ i\in\{1,\ldots,p\},\\
    4(p-q+1)e_i &\text{if} \ \ i\in\{p+1,\ldots,p+q\}.
  \end{cases}
\end{align*}

The linear operator $\DX$ is invariant under orthogonal
transformations of the elements in $\X$ in the following sense:

\begin{lemma}
  Let $\X=\{x_1,x_2,\ldots,x_d\}$ and $\X'=\{x_1',\ldots,x_d'\}$ be
  subsets of a Lie algebra $\g$ such that $x_i'=\sumjd R_{ij}x_j$ for
  some orthogonal $d\times d$-matrix $R$. Then $\DX(a)=\Delta_{\X'}(a)$ for all
  $a\in\g$.
\end{lemma}

\begin{proof}
  The proof is given by the following calculation:
  \begin{gather*}
      \Delta_{\X'}(a)  = \sum_{j,k,l=1}^dR_{jk}R_{jl}\Ccom{a}{x_k}{x_l}
      = \sum_{k,l=1}^d\paraa{R^TR}_{kl}\Ccom{a}{x_k}{x_l}\\
\hphantom{\Delta_{\X'}(a)}{} = \sum_{k,l=1}^d\delta_{kl}\Ccom{a}{x_k}{x_l}
      = \sumjd\Ccom{a}{x_j}{x_j} = \DX(a).\tag*{\qed}
    \end{gather*}\renewcommand{\qed}{}
\end{proof}

\begin{remark}
  Note that it is not necessarily true that $(\g,\X')$ is a \DMSA\
  when $(\g,\X)$ is a \DMSA. However, if we let
  $\Spec\paraa{(\g,\X)}=\{\mu_1,\ldots,\mu_k\}$ and
  $m_1,\ldots,m_k$ be the multiplicities of each eigenvalue, then any
  block-diagonal orthogonal matrix $R=\diag(R_1,\ldots,R_k)$ (where
  $R_i$ has dimension $m_i$) will generate a \DMSA. In other words, we
  can always choose to make an orthogonal transformation among those
  $x_i$ that belong to the same eigenvalue.
\end{remark}

  The preceding lemma enables us to make the following
observation. Let $\X$ be a set of linearly independent elements in a
$n$-dimensional Lie algebra $\g$, and let $\vspan{\X}$ denote the
linear span of the elements in $\X$. Furthermore, assume that
$\vspan{\X}$ is closed under the action of $\DX$,
i.e.\ $\DX(x)\in\vspan{\X}$ for all $x\in\vspan{\X}$. Relative to a
basis where $x_1,\ldots,x_d$ are chosen to be the f\/irst $d$ basis
elements, the $n\times n$ matrix of $\DX$ has the block form
\begin{gather*}
  \DX =
  \begin{pmatrix}
    X_0 & A \\
    0 & B
  \end{pmatrix}
\end{gather*}
where $X_0$ is a $d\times d$ matrix. If $X_0$ is diagonalizable by an
orthogonal matrix $R$ then the elements $x_i'=\sum_{j=1}^dR_{ij}x_j$ will be
eigenvectors of $\Delta_{\X'}$, since the action of $\DX$ is invariant
under orthogonal transformations in $\X$. Thus, $(\g,\X')$ is a
\DMSA. In particular, if we choose an orthonormal basis of $\g$, then
the matrices $\ad_{x_i}$ are antisymmetric, which implies that
$\DX$ is symmetric. In this case, $X_0$ will be diagonalizable by an
orthogonal matrix.

We will now concentrate on subsets $\X$ of the Lie algebra $\sl_n$,
such that $\A=(\sl_n,\X)$ is a \DMSA. To perform calculations the
following set of conventions will be used:
$\alpha_1,\ldots,\alpha_{n-1}$ denote the simple roots and for every
positive root $\alpha$, we choose elements
$e_{\alpha}$, $e_{-\alpha}$, $h_{\alpha}$ such that
\begin{gather*}
 [h,e_\alpha]=\alpha(h)e_\alpha,\qquad
 [e_{\alpha},e_{-\alpha}]=h_{\alpha},
\end{gather*}
and $h_\alpha$ is the element of the Cartan subalgebra $\h$ such
that $\alpha(h)=K(h_\alpha,h)$ for all $h\in\h$.  For any pair of
roots $\alpha,\beta$ we def\/ine the constants $N(\alpha,\beta)$ by
\begin{gather*}
  [e_\alpha,e_\beta]=N(\alpha,\beta)e_{\alpha+\beta},
\end{gather*}
and when $\alpha+\beta$ is not a root, we set $N(\alpha,\beta)=0$. In
$\sl_n$ all roots have the same length, and we denote
$(\alpha,\alpha)\equiv K(\ha,\ha)=\alpha(h_\alpha)=l^2$. With these conventions, the
constants $N(\alpha,\beta)$ satisfy the relations
\begin{gather*}
   N(\alpha,\beta)=N(\beta,\gamma)=N(\gamma,\alpha)\qquad\text{if} \ \ \alpha+\beta+\gamma=0,\\
   N(\alpha,\beta)N(-\alpha,-\beta)=-\frac{l^2}{2}q(p+1),
\end{gather*}
where $p$, $q$ are positive integers such that
$\beta-p\alpha,\ldots,\beta,\ldots,\beta+q\alpha$ are
roots. Furthermore, in~$\sl_n$ it holds that if $\beta+\alpha$ is a
root then $\beta-\alpha$ is not a root and $\beta\pm 2\alpha$ is never a
root. Therefore, if~$N(\alpha,\beta)$ is non-zero, we have that
\begin{gather*}
  N(\alpha,\beta)N(-\alpha,-\beta)=-\frac{1}{2} l^2.
\end{gather*}
Although the following result does not depend on it, we will for
def\/initeness choose each $N(\alpha,\beta)$ such that
$N(-\alpha,-\beta)=-N(\alpha,\beta)$.

\begin{lemma}\label{lemma:ealphaplus}
  For every positive root $\alpha$ in $\sl_n$, we set
  \begin{gather*}
    \eap = ic\paraa{\ea+\ema}\qquad\text{and}\qquad
    \eam = c\paraa{\ea-\ema},
  \end{gather*}
  for an arbitrary $c\in\reals$. Then the following holds
  \begin{enumerate}\itemsep=0pt
  \item $\Ccom{\eap}{\ebp}{\ebp} = \Ccom{\eap}{\ebm}{\ebm}=-\frac{1}{2} c^2l^2\eap$
    \qquad\text{$($when $\alpha\pm\beta$ is a root$)$},
  \item $\Ccom{\eam}{\ebp}{\ebp} = \Ccom{\eam}{\ebm}{\ebm}=-\frac{1}{2} c^2l^2\eam$
    \qquad\text{$($when $\alpha\pm\beta$ is a root$)$},
  \item $\Ccom{e^{\pm}_\alpha}{e^{\mp}_\alpha}{e^{\mp}_\alpha} = -2 c^2l^2e^{\pm}_{\alpha}$,
  \item $\Ccom{e^{\pm}_\alpha}{\hb}{\hb}=(\alpha,\beta)^2e^{\pm}_\alpha$,
  \item $\Ccom{\ha}{e^{\pm}_\beta}{e^{\pm}_\beta} = \mp 2 c^2(\alpha,\beta)\hb$.
  \end{enumerate}
\end{lemma}

From this lemma, it is easy to construct a couple of examples.

\begin{example}
  Let $\X=\{e^{\pm}_{\beta_1},\ldots,e^{\pm}_{\beta_d}\}$ (where the
  signs are chosen independently) for any positive roots $\beta_i$. In
  this case, $\Ccom{x_i}{x_j}{x_j}$ is proportional to $x_i$ for all
  $x_i,x_j\in\X$.
\end{example}

\begin{example}
  Let
  $\X=\{h_{\beta_1},\ldots,h_{\beta_k},\ep_{\gamma_1},\emi_{\gamma_1},\ldots,\ep_{\gamma_l},\emi_{\gamma_l}\}$.
  Now, $\Ccom{h_{\beta_i}}{\ep_{\gamma_j}}{\ep_{\gamma_j}}$ might not
  be proportional to $h_{\beta_i}$. However, since both
  $\ep_{\gamma_j},\emi_{\gamma_j}\in\X$ this term will cancel against
  $\Ccom{h_{\beta_i}}{\emi_{\gamma_j}}{\emi_{\gamma_j}}$. Thus, $\DX(h_{\beta_i})=0$ for $i=1,\ldots,k$.
\end{example}

\begin{example}
  Let
  $\X=\{h_{\beta_1},\ldots,h_{\beta_k},e^{\pm}_{\beta_1},\ldots,e^{\pm}_{\beta_k}\}$
  (where the signs are chosen to be the same). In this case
  $\DX(h_{\beta_i})$ will not be proportional to
  $h_{\beta_i}$. However, the matrix $\DX$ will be symmetric, which
  implies that there exists an orthogonal $k\times k$ matrix $R$ such
  that
  $(\sl_n,\{x_1,\ldots,x_k,e^{\pm}_{\beta_1},\ldots$, $e^{\pm}_{\beta_k}\})$
  is a \DMSA\ if $x_i=\sum_{j=1}^kR_{ij}h_{\beta_j}$.
\end{example}

When the dimension $d=2m$ (i.e.\ even) and every eigenvalue in $\Spec(\A)$ has an
even multiplicity (which is relevant for one of the applications
mentioned in the introduction) there is a~convenient complexif\/ied
basis provided by
\begin{gather*}
   t_i = x_{2i-1}+ix_{2i},\qquad
   s_i = x_{2i-1}-ix_{2i}
\end{gather*}
for $i=1,\ldots,m$. The def\/ining relations of a $\DMSA$ can then be written as
\begin{gather*}
   2\mu_it_i = \sum_{j=1}^m\big(\Ccom{t_i}{s_j}{t_j}+\Ccom{t_i}{t_j}{s_j}\big),\qquad
  2\mu_is_i = \sum_{j=1}^m\big(\Ccom{s_i}{t_j}{s_j}+\Ccom{s_i}{s_j}{t_j}\big).
\end{gather*}

 The lowest dimensional non-trivial \DMSA\ has dimension 2. In this case
the algebra is gene\-ra\-ted by $x$ and $y$ satisfying
\begin{gather*}
   \Ccom{x}{y}{y}=\lambda x,\qquad
   \Ccom{y}{x}{x}=\mu y,
\end{gather*}
and by def\/ining $z=-i[x,y]$ one sees that $x$, $y$ and $z$ span a 3-dimensional
Lie algebra. By rescaling the elements we obtain the following result:
\begin{itemize}\itemsep=0pt
\item $\lambda\neq 0,\mu\neq 0$: $\A$ is isomorphic to $\sl_2$,
\item $\lambda=\mu=0$: $\A$ is isomorphic to the Heisenberg algebra,
\item $\lambda\neq 0,\mu=0$ or $\lambda=0,\mu\neq 0$: $\A$ is
  isomorphic to the Lie algebra $\text{VII}_1$ in the Bianchi
classif\/ication \cite{b:threeDimLieAlgebras}. This algebra is def\/ined
  by the relations: $[u,v]=-w$, $[v,w]=0$ and $[w,u]=-v$.
\end{itemize}

\section[Hermitian representations of $\U_4(\A)$]{Hermitian representations of $\boldsymbol{\U_4(\A)}$}\label{sec:hermitianReps}

 In general, any representation of the Lie algebra $\g$ gives
rise to a representation of the $\DMSA$ $(\g,\X)$. In the following,
we shall however concentrate on f\/inding \emph{Hermitian}
representations.  Hermitian representations $\phi$ of $\U_d(\A)$ are
given by $\phi(x_i)=X_i$ where $X_1,\ldots,X_d$ are Hermitian matrices
satisfying
\begin{gather*}
  \sumjd\Ccom{X_i}{X_j}{X_j}=\mu_iX_i\qquad\text{for} \ \ i=1,\ldots,d.
\end{gather*}
We note that, unless all $\mu_i$ are real, no Hermitian (or
anti-Hermitian) representations can exist (except for the trivial one:
$\phi(x_i)=0$ for all $i$). Hence, from now on we will assume the
spectrum to be real and all representations to be Hermitian.

Since the subalgebra (of the full matrix algebra) generated by the
matrices $\{\phi(x_1),\ldots,\phi(x_d)\}$ (through arbitrary products
and sums) is invariant under Hermitian conjugation, the following
result is immediate.

\begin{proposition}\label{prop:completelyReducible}
  Any Hermitian representation of $\U_d(\A)$ is completely reducible.
\end{proposition}

 As \DMSA s of dimension 2 are isomorphic to Lie algebras, we
will start by considering the case when $d=4$. We expect that these
algebras have a rich structure of representations even for the case
when $\mu_1=\cdots=\mu_4$. Namely, since the equations def\/ining a
$4$-dimensional \DMSA\ can be thought of as discrete analogues of
minimal surface equations in $S^3$ (see the introduction), and minimal
surfaces of any genus exist in $S^3$ \cite{l:completeMinimalSurfaces},
we believe that there will be representations corresponding to many
(if not all) of these surfaces.

Since the def\/ining relations of $\U_d(\A)$ are expressed entirely in
terms of commutators, the tensor product of Lie algebra
representations, i.e.
\begin{gather*}
  \paraa{\phi\otimes\phi'}(x) = \phi(x)\otimes\mid + \mid\otimes\phi'(x),
\end{gather*}
also def\/ines a tensor product for representations of \DMSA s.  In
contrast to Lie algebras, the tensor product of two irreducible
representations might again be irreducible (as we shall explicitly see
for $\U_4(\A)$).  Thus, when studying the representation theory of
$\U_d(\A)$ it becomes natural to look, not only for irreducible
representations, but also for \emph{prime representations},
i.e.\ irreducible representations that can not be written as a tensor product of
two other representations.

To each Hermitian representation $\phi$ of $\U_4(\A)$ we shall
associate a directed graph with vertices in $\complex$, such that the
vertices of the graph are placed at the characteristic roots of
$\phi(x_1+ix_2)$. The edges of the graph are determined by the matrix
$\phi(x_3+ix_4)$ as described below. We note that this construction
can be carried out for Hermitian representations of any algebra on at
most four generators. In the following, we will use the notation
$t_1=x_1+ix_2$, $t_2=x_3+ix_4$, $\phi(t_1)=\Lambda$ and
$\phi(t_2)=T$. Let us start by recalling the directed graph of a
matrix.

\begin{definition}
  Let $T$ be a $n\times n$ matrix and let $G=(V,E)$ be a directed
  graph on $n$ vertices with vertex set $V=\{1,\ldots,n\}$ and edge
  set $E\subseteq V\times V$. We say that $G$ is the \emph{directed
    graph of~$T$}, and write $G=G_T$, if it holds that
  \begin{align*}
    T_{ij}\neq 0 \equivalent (i,j)\in E.
  \end{align*}
  for $i,j=1,\ldots,n$.
\end{definition}

 The idea is now to associate a graph to every
representation, such that each vertex is assigned an eigenvalue of
$\Lambda$ and the graph itself being the directed graph of~$T$. Needless to say, the graph of~$T$ depends on the basis chosen and
therefore we will introduce a particular choice of basis in which all
graphs will be refered to.

\begin{definition}
  Let $\Lambda$ and $T$ be linear operators on $W=\complex^n$, and let
  $\B$ be the $\ast$-algebra generated by
  $\Lambda$, $\Lambda^\dagger$, $T$, $\Td$. Furthermore, let $W=W_1\oplus
  \cdots\oplus W_m$ be a decomposition of $W$ into irreducible
  subspaces with respect to $\B$. For each $i$, let
  $v_1^{(i)},\ldots,v_{n_i}^{(i)}$ denote a Jordan basis for
  $\Lambda\big|_{W_i}$. Then
  $v_1^{(1)},\ldots,v_{n_1}^{(1)},\ldots,v_1^{(m)},\ldots,v_{n_m}^{(m)}$
  is a basis for $W$ and is called \emph{a Jordan basis of $\Lambda$
    with respect to~$T$}.
\end{definition}

\begin{definition}
  Let $\phi$ be a $n$-dimensional representation of
  $\U_4(\A)$ and let $v_1,\ldots,v_n$ denote a~Jordan basis of
  $\Lambda=\phi(t_1)$ with respect to $T=\phi(t_2)$.  Def\/ine
  the matrix $\alpha$ by def\/ining its matrix elements through
  \begin{align*}
    Tv_i=\sum_{j=1}^n\alpha_{ji}v_j,
  \end{align*}
  i.e., $\alpha$ is the matrix of $T$ in the given Jordan basis.
  Furthermore, let $G_\alpha=(\{1,2,\ldots,n\},E)$ denote the directed
  graph of $\alpha$ and let $\lambda:V\to\complex$ be def\/ined by
  $\lambda(i)=\lambda_i$, where $\lambda_i$ is the eigenvalue
  corresponding to $v_i$. We set $G_\phi=(G_\alpha,\lambda)$ and call
  $G_\phi$ \emph{a representation graph of $\phi$}.
\end{definition}

 Two representation graphs $G_\phi=(\{1,\ldots,n\},E,\lambda)$ and
$G_{\phi'}=(\{1,\ldots,n\},E',\lambda')$ are \emph{isomorphic} if
there exists a permutation $\sigma\in S_n$ such that $(i,j)\in
E\Leftrightarrow \paraa{\sigma(i),\sigma(j)}\in E'$ and
$\lambda=\lambda'\circ\sigma$. In particular, $(\{1,\ldots,n\},E)$ and
$(\{1,\ldots,n\},E')$ are isomorphic as directed graphs.

Note that two representation graphs corresponding to the
same representation need not be isomorphic. This could be resolved by
further f\/ixing the basis in which the directed graph of $T$ is
calculated. However, let us postpone this choice and study the
properties of representation graphs that follow from the above
def\/inition.

The f\/irst property that one might wish for, is a correspondence
between disconnected components of the representation graph and the
irreducible components of the representation. That a connected graph
corresponds to an irreducible representation follows immediately from
the def\/inition of the Jordan basis with respect to $T$.

\begin{proposition}\label{prop:connectedIrreducible}
  Let $G_\phi$ be a representation graph of $\phi$. If $G_\phi$ is
  connected then $\phi$ is irreducible.
\end{proposition}

\begin{proof}
  Let $G_\phi$ be a connected representation graph of $\phi$. If
  $\phi$ is reducible, then $G_\phi$ consists of at least two
  components, since the matrix $\alpha$ is block diagonal with at
  least two blocks, by the construction of the Jordan basis. Hence,
  $\phi$ must be irreducible.
\end{proof}

 For convenience, let us introduce some terminology
indicating when the matrices of a representation have certain
properties.

\begin{definition}
  Let $\phi$ be a representation of $\U_4(\A)$. If
  $\Lambda=\phi(t_1)$ is diagonalizable then $\phi$ is called
  diagonalizable. If all eigenvalues of $\Lambda$ are distinct, then
  $\phi$ is called non-degenerate. If $\Lambda$ is normal then $\phi$ is
  called semi-normal. If both $\Lambda$ and $T=\phi(t_2)$ are
  normal then $\phi$ is called normal. If $\Lambda$ and $T$ are
  unitary, then $\phi$ is called unitary.
\end{definition}

 For semi-normal representations, the result in Proposition
\ref{prop:connectedIrreducible} can be strengthened to an if and
only if statement.

\begin{proposition}\label{prop:semiNormalConnected}
  Let $\phi$ be a semi-normal representation of $\U_4(\A)$ and let
  $G_\phi$ be a representation graph of $\phi$. Then $\phi$ is
  irreducible if and only if $G_\phi$ is connected.
\end{proposition}

\begin{proof}
  Assuming that $G_\phi$ is connected, the f\/irst implication follows
  from Proposition \ref{prop:connectedIrreducible}. Now, assume that
  $\phi$ is an irreducible $n$-dimensional representation. When
  $\Lambda$ is normal, the Jordan basis is given by a set of
  orthogonal vectors $\{v_1,\ldots,v_n\}$ (the eigenvectors of
  $\Lambda$). Therefore, any matrix $P$, bringing $\Lambda$ to the
  (diagonal) Jordan normal form $P^{-1}\Lambda P$, can be written as
  $P=UD$ where $U$ is a unitary matrix and $D$ is an invertible
  diagonal matrix (ref\/lecting the choice of length of the
  eigenvectors). In such a basis, the matrix $\alpha$, representing
  $T$, and the matrix $\alphat$, representing $T^\dagger$ are related
  by $\alphat^\dagger=D^\dagger D\alpha(D^\dagger D)^{-1}$. Clearly,
  since conjugation by a diagonal invertible matrix does not change
  the structure of non-zero matrix elements, the graph of $\alphat$ is
  obtained from the graph of~$\alpha$ by reversing the arrows. In
  particular, this implies that~$G_{\alphat}$ and~$G_{\alpha}$ have
  the same number of connected components. We now continue to prove
  that if $G_\alpha$ is disconnected then there exists an invariant
  subspace, which contradicts the assumption that~$\phi$ is
  irreducible.

  Assume that $G_\alpha$ is disconnected and let
  $\{i_1,\ldots,i_k\}\subset V$ (with $k<n$) be the vertices
  corresponding to one of the components. By the construction of
  $G_\alpha$, the vectors $v_{i_1},\ldots,v_{i_k}$ span an invariant
  subspace for $\alpha$. Since $G_{\alphat}$ is obtained from
  $G_\alpha$ by reversing all arrows, this is also an invariant
  subspace for $\alphat$. Moreover, since $v_{i_1},\ldots,v_{i_k}$ are
  eigenvectors of both $\Lambda$ and $\Lambda^{\dagger}$, the space
  spanned by these vectors is also an invariant subspace of $\Lambda$
  and $\Lambda^\dagger$. This implies that the representation is
  reducible, which contradicts the assumption. Hence, $G_\phi$ is
  connected.
\end{proof}

 When the eigenvalues of $\Lambda=\phi(t_1)$ are distinct, one can
easily show that all representation graphs of $\phi$ are isomorphic.

\begin{proposition}\label{prop:distinctEigenvaluesBasis}
  Let $\phi$ be a non-degenerate representation. Then any Jordan basis for
  $\Lambda$ is a~Jordan basis for $\Lambda$ with respect to $T$ up to
  a permutation of the basis vectors. Moreover, all representation
  graphs of $\phi$ are isomorphic.
\end{proposition}

\begin{proof}
  When all eigenvalues of $\Lambda$ are distinct, the only freedom in
  choosing a basis in which~$\Lambda$ is diagonal lies in the length
  of the eigenvectors and the ordering of the basis vectors. Hence,
  given two Jordan bases for $\Lambda$ it is always possible to apply
  a permutation to obtain one basis from the other, up to a rescaling
  of the vectors. Furthermore, a rescaling of the basis vectors does
  not change the block diagonal form of a matrix. Hence, any Jordan
  basis of $\Lambda$ is a Jordan basis of $\Lambda$ with respect to~$T$ up to a permutation. In particular, this implies that any two
  representation graphs are related by a permutation of the vertices.
\end{proof}

 Proposition \ref{prop:distinctEigenvaluesBasis} has the
consequence that if one constructs a representation $\phi$, in which
$\Lambda$ is diagonal and has distinct eigenvalues, then the directed
graph of $T$ is the unique representation graph of $\phi$.

Let us now study the representation graph of the tensor product.
For directed graphs, forming the tensor product
\begin{gather*}
  \hat{T} = T\otimes\mid + \mid\otimes T'
\end{gather*}
amounts to taking the \emph{Cartesian product} of $G_T$ and $G_{T'}$
\cite{s:graphMultiplication,ht:connectednessGraph}. The Cartesian
product of two graphs $G=(V,E)$ and $H=(U,F)$ is def\/ined as the graph
$G'=(V\times U,E')$ such that
\begin{align*}
  &\paraa{(v_1,u_1),(v_2,u_2)}\in E'\equivalent\\
  &\big\{v_1=v_2\text{ and }(u_1,u_2)\in F\big\}
  \text{ or }
  \big\{u_1=u_2\text{ and }(v_1,v_2)\in E\big\}.
\end{align*}
Now, one might ask if the Cartesian product of two representation
graphs is a representation graph of the tensor product? This is not
always true, but we have the following result.

\begin{proposition}\label{prop:tensorProductGraph}
  Let $\phi$ and $\phi'$ be representations such that
  $\phi\otimes\phi'$ is a non-degenerate representation. Then $\phi$
  and $\phi'$ are non-degenerate and $G_{\phi\otimes\phi'}$ is the
  Cartesian product of~$G_\phi$ and~$G_{\phi'}$.
\end{proposition}

\begin{proof}
  Let $\L=\phi(t_1)$ and $\L'=\phi'(t_1)$ and let $P$ and $Q$ be
  matrices whose column vectors are Jordan bases of $\L$ and $\L'$
  with respect to $T=\phi(t_2)$ and $T'=\phi'(t_2)$. By assumption,
  the eigenvalues of $\hat{\L}=\L\otimes\mid + \mid\otimes\L'$ are
  distinct, which implies that the eigenvalues of the matrix
  \begin{gather*}
    M = \paraa{P\otimes Q}^{-1}\bracketb{\L\otimes\mid + \mid\otimes\L'}\paraa{P\otimes Q}
    =\paraa{P^{-1}\L P}\otimes\mid + \mid\otimes\paraa{Q^{-1}\L'Q}
  \end{gather*}
  are distinct. Since $P^{-1}\L P$ and $Q^{-1}\L'Q$ are upper
  triangular (and hence has their eigenvalues on the diagonal) the
  matrix $M$ will also be upper triangular. The diagonal elements of
  $M$ (its eigenvalues) will be all possible sums of eigenvalues from
  $\L$ and $\L'$. Since the eigenvalues of~$M$ are distinct the
  eigenvalues of $\L$ and $\L'$ must be distinct. This proves the
  f\/irst part of the statement.

  Since $\L$ and $\L'$ have distinct eigenvalues, the matrices
  $P^{-1}\L P$ and $Q^{-1}\L'Q$ are in fact diagonal, which implies
  that the matrix $M$ is diagonal, so $P\otimes Q$ clearly provides us
  with a~Jordan basis for $\hat{\L}$. Since $\phi\otimes\phi'$ is
  non-degenerate, Proposition~\ref{prop:distinctEigenvaluesBasis} tells
  us that the (unique) representation graph is given by the directed graph of
  \begin{gather*}
    \paraa{P\otimes Q}^{-1}\big[T\otimes\mid + \mid\otimes T'\big]\paraa{P\otimes Q}
    =\paraa{P^{-1}T P}\otimes\mid + \mid\otimes\paraa{Q^{-1}T'Q}.
  \end{gather*}
  Now, since $P^{-1}T P$ and $Q^{-1}T'Q$ def\/ine $G_\phi$ and
  $G_{\phi'}$, we conclude that $G_{\phi\otimes\phi'}$ is given by the
  Cartesian product of $G_\phi$ and $G_{\phi'}$.
\end{proof}

 Let us now proceed to construct representations of
$\U_4(\A)$. As noted earlier, even the case when
$\mu_1=\cdots=\mu_4\neq 0$ is expected to have a rich representation
theory. Therefore, we will start by concentrating on the case for
which $\mu_1=\mu_2=\mu$ and $\mu_3=\mu_4=\rho$
(which is also relevant for applications, as mentioned in the
introduction). In this case, representations are found by solving the
matrix equations
\begin{gather*}
   2\mu\L = \Ccom{\L}{T}{\Td} + \Ccom{\L}{\Td}{T} + \Ccom{\L}{\Ld}{\L},\\
   2\rho T = \Ccom{T}{\L}{\Ld} + \Ccom{T}{\Ld}{\L} + \Ccom{T}{\Td}{T}.
\end{gather*}
The action of the group $O(2)\times O(2)$ can be explicitly realized
by letting $\Lambda\to e^{i\theta}\Lambda$ and $T\to e^{i\theta'}T$,
which gives a new representation for any $\theta,\theta'\in\reals$;
this representation will be denoted by $\phitt$ and is in general not
equivalent to $\phi$ since the eigenvalues of e.g.~$\Lambda$ will be
dif\/ferent.  This enables us to construct new irreducible
representations from a given one via the tensor product. Namely, let
$\phi$ be a non-degenerate irreducible representation; then one can
always choose $\theta$, $\theta'$ such that $\phi\otimes\phitt$ is a
non-degenerate representation. By Proposition~\ref{prop:tensorProductGraph} the representation graph of
$\phi\otimes\phitt$ will be the Cartesian product of two connected
graphs (the representation graphs of~$\phi$ and~$\phitt$), which
implies that it is connected~\cite{ht:connectednessGraph}. Hence, it
follows from Proposition~\ref{prop:connectedIrreducible} that
$\phi\otimes\phitt$ is irreducible.

\subsection{The Fuzzy sphere}

  As any semi-simple Lie algebra is itself a \DMSA, it follows
that Hermitian representations of $\su(2)$ should induce Hermitian
representations of $\U_4(\A)$. Indeed, choosing Hermitian $n\times n$
matrices $S_1$, $S_2$, $S_3$, with non-zero elements
\begin{gather*}
   \paraa{S_1}_{k,k+1}=\frac{1}{2}\sqrt{k(n-k)}=\paraa{S_1}_{k+1,k},\qquad k=1,\ldots,n-1,\\
   \paraa{S_2}_{k,k+1}=-\frac{i}{2}\sqrt{k(n-k)}=-\paraa{S_2}_{k+1,k},\qquad k=1,\ldots,n-1,\\
   \paraa{S_3}_{k,k} = \frac{1}{2}(n+1-2k),\qquad k=1,\ldots,n,
\end{gather*}
satisfying $[S_i,S_j]=i\epsilon_{ijk}S_k$, yields a representation
$\phi$ by def\/ining
\begin{gather*}
   \L =\phi(t_1) = e^{i\theta}S_3,\qquad
   T = \phi(t_2) = S_1+iS_2.
\end{gather*}
for any $\theta\in\reals$ (no $\theta'$ appears in $\phi(t_2)$ since it
can always be removed by conjugating with a diagonal unitary
matrix). One easily calculates that this is a representation of
$\U_4(\A)$ with $\Spec(\A)=\{2\}$.  Note that in the special case when
$\theta=0$, in which case $\L$ is Hermitian, this provides a
representation of $\U_3(\A)$.

\begin{figure}[h]
  \centering
  \includegraphics[width=8cm]{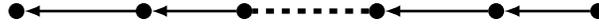}
  \caption{The representation graph of the Fuzzy sphere.}
  \label{fig:fuzzyspheregraph}
\end{figure}

This is a non-degenerate semi-normal representation, and the
representation graph takes the form as in Fig.~\ref{fig:fuzzyspheregraph}.

 Furthermore, this representation is irreducible by
Proposition \ref{prop:connectedIrreducible}, and its representation
graph is prime with respect to the Cartesian product since any
Cartesian product graph with $n$ vertices has at least $n$ edges
(whereas the above graph has $n-1$ edges). The matrix algebras
generated by these matrices have been used to construct sequences of
matrix algebras (of increasing dimension) converging to the Poisson
algebra of smooth functions on $S^2$ \cite{h:phdthesis}.

For increasing $n$, the algebras that
are generated by these matrices (with an appropriate normalization)
are recognized as a sequence converging to the Poisson algebra of
functions on~$S^2$~\cite{h:phdthesis}.

\begin{figure}[h]
  \centering
  \includegraphics[width=8cm]{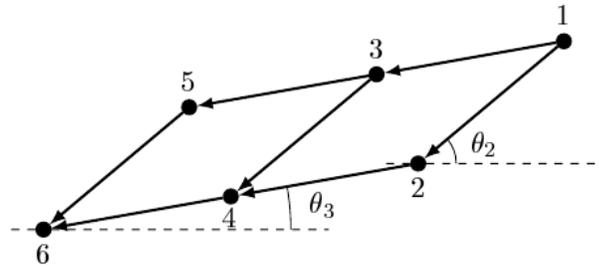}
  \caption{The representation graph of a tensor product of two Fuzzy sphere representations.}
  \label{fig:spheretensorgraph}
\end{figure}

Let us for this case demonstrate the tensor product and construct the
corresponding Cartesian product of the representation graphs. For
simplicity, we let $\phi_2$ and $\phi_3$ be a two- respectively
three-dimensional representation of the type described above, and set
\begin{gather*}
   \phi(t_1) = \phi_3(t_1)\otimes\mid_2 + \mid_3\otimes\phi_2(t_1), \\
   \phi(t_2) = \phi_3(t_2)\otimes\mid_2 + \mid_3\otimes\phi_2(t_2).
\end{gather*}
If we denote the arbitrary phases by $\theta_2$ and $\theta_3$ the
representation graph takes the form as in Fig.~\ref{fig:spheretensorgraph}.

 Since this representation is non-degenerate it will be irreducible by
Proposition~\ref{prop:connectedIrreducible}, and in general we obtain
inequivalent representations for dif\/ferent choices of $\theta_2$ and
$\theta_3$.

\begin{figure}[h]
  \centering
  \includegraphics[width=10cm]{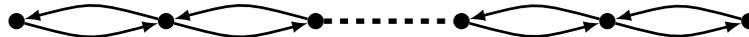}
  \caption{The representation graph of a normal representation constructed from
    $\su(2)$.}
\label{fig:normalfuzzyspheregraph}
\end{figure}

We note that the matrices $S_1$, $S_2$, $S_3$ gives rise to another
representation by setting
\begin{gather}
   \L = zS_3,\label{eq:sphereNormalRep1}\\
   T = w(S_1+aS_2)\label{eq:sphereNormalRep2}
\end{gather}
for arbitrary $z,w\in\complex$ and $a\in\reals$. This gives a
representation of $\U_4(\A)$ with
$\Spec(\A)=\{|z|^2,|w|^2(1+a^2)\}$. The representation graph can be seen in
Fig.~\ref{fig:normalfuzzyspheregraph}.

 We conclude that this is an irreducible non-degenerate
normal representation, which is not equivalent to the Fuzzy sphere,
since the corresponding graphs are not isomorphic. Moreover, its
representation graph is prime with respect to the Cartesian product.

\subsection{The Fuzzy torus}

\begin{figure}[h]
  \centering
  \includegraphics[width=4cm]{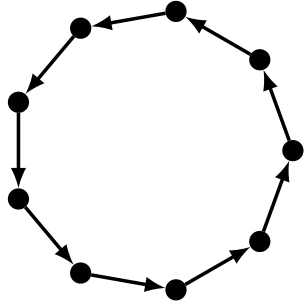}
  \caption{The representation graph of the Fuzzy torus.}
  \label{fig:fuzzytorusgraph}
\end{figure}

The fuzzy torus algebra
(cp.~\cite{ffz:trigonometric,h:diffeomorphism}) is generated by the
matrices $g$ and $h$, with non-zero elements
\begin{gather*}
   (h)_{k,k+1} = 1,\qquad (h)_{n,1}=1,  \qquad k =1,\ldots,n-1,\\
   (g)_{kk} = q^{k-1},\qquad  k=1,\ldots,n,
\end{gather*}
fulf\/illing the relation $hg=q\cdot gh$ with $q^n=1$. It is known that they generate
matrix sequences that converge to functions on $T^2$. A representation $\phi$ of $\U_4(\A)$, with
$\Spec(\A)=\{|1-q|^2/2\}$, is obtained by setting
\begin{gather*}
   \Lambda=\phi(t_1) = e^{i\theta}g, \qquad
   T=\phi(t_2) = e^{i\theta'}h,
\end{gather*}
for any $\theta,\theta'\in\reals$.  This is an irreducible
non-degenerate unitary representation, with a representation graph as
in Fig.~\ref{fig:fuzzytorusgraph}.

 Furthermore, this graph is prime with respect to the
Cartesian product. Let us now show that this is essentially the only
irreducible unitary representation when $\Spec(\A)=\{\mu\}$.

\subsection{Unitary representations}

When $\Lambda$ and $T$ are unitary and $\Spec(\A)=\{\mu\}$, the
equations can be written as
\begin{gather}
   \lambda\Lambda = \Td\L T + T\L\Td,\label{eq:unitaryEq1}\\
   \lambda T = \Ld T\L + \L T\Ld,\label{eq:unitaryEq2}
\end{gather}
where we have introduced $\lambda=2-\mu$. By multiplying the f\/irst
equation from the right by $T$ and the second equation from the left
by $\L$ we note that $[\L T,T\L]=0$. Hence, by a unitary
transformation, we can always choose a basis such that $D=\L T$ and
$\Dt=T\L$ are diagonal. Let us denote the eigenvalues by
\begin{gather*}
   D = \diag(d_1,\ldots,d_n)=\diag\paraa{e^{i\vphi_1},\ldots,e^{i\vphi_n}},\\
   \Dt = \diag(\dt_1,\ldots,\dt_n)=\diag\paraa{e^{i\vphit_1},\ldots,e^{i\vphit_n}}.
\end{gather*}
The equations (\ref{eq:unitaryEq1}) and (\ref{eq:unitaryEq2})
(together with $D=\L T$ and $\Dt=T\L$) can equivalently be written~as
\begin{gather}
  \lambda\L = \L\Dtd D+\Dt\L\Dtd,\label{eq:unitaryDDt1}\\
  \lambda\L\Dtd=\Dtd\L + \L D^\dagger,\label{eq:unitaryDDt2}\\
  \L\Dt=D\L,\label{eq:unitaryDDt3}\\
  T=\Dt\Ld.\label{eq:unitaryDDt4}
\end{gather}
Thus, given unitary $\Lambda$, $D$, $\Dt$ satisfying
(\ref{eq:unitaryDDt1})--(\ref{eq:unitaryDDt3}), we obtain a solution
to the original equations by def\/ining $T=\Dt\Ld$. Written out in
components, the three f\/irst equations become
\begin{gather}
 \L_{ij}\bracketb{\lambda-\bar{\dt}_j d_j-\dt_i\bar{\dt}_j}=0,\label{eq:unitaryComp1}\\
 \bar{\L}_{ij}\bracketb{\lambda\dt_j-\dt_i-d_j}=0,\label{eq:unitaryComp2}\\
\L_{ij}\bracketa{\dt_j-d_i}=0.\label{eq:unitaryComp3}
\end{gather}
If $\L_{ij}\neq 0$ then we obtain the following relations
\begin{gather*}
   \twovec{d_j}{\dt_j}=\twomatrix{\lambda}{-1}{1}{0}\twovec{d_i}{\dt_i}\equiv s\twovec{d_i}{\dt_i}
\end{gather*}
since (\ref{eq:unitaryComp2}) and (\ref{eq:unitaryComp3}) together
imply (\ref{eq:unitaryComp1}). Now, consider the directed graph
$G_\L=(V,E)$ of $\L$, where we have assigned the vector
$\xv_i=(d_i,\dt_i)$ to each vertex $i\in V$. We can restrict ourselves
to connected graphs, since if $G_\L$ is disconnected then the
representation will trivially be reducible. The above considerations
tell us that whenever there is an edge $(i,j)\in E$, it must hold that
$\xv_j=s(\xv_i)$. In particular, since $D$ and $\Dt$ are unitary
matrices, the map $s$ must take $\xv\in S^1\times S^1$ to another
vector in $S^1\times S^1$. If $\xv=(e^{i\vphi},e^{i\vphit})$ then this
is true only if $\lambda=0$ or
\begin{gather}
  \lambda=2\cos(\vphi-\vphit).\label{eq:philambda}
\end{gather}
This observation leads to the following result.

\begin{proposition}
  Let $\A$ be a \DMSA\ with $\Spec(\A)=\{\mu\}$. If $\mu<0$ or
  $\mu>4$ then there exists no unitary representation of $\U_4(\A)$.
\end{proposition}

\begin{proof}
  Assume that $\Lambda$ and $T$ provides a unitary representation and
  that a basis has been chosen in which $D$ and $\Dt$ are diagonal. Since
  $\Lambda$ is a unitary matrix at least one of its matrix element has to be
  non-zero, say $\L_{ij}\neq 0$. Then equation (\ref{eq:philambda})
  must hold, which is impossible if $\lambda<-2$ or $\lambda>2$.
\end{proof}

 From the above result it follows that whenever a unitary
representation exists, then there exists a $\beta\in[0,\pi/2]$ such
that $\lambda=2\cos 2\beta$. We will now proceed in analogy with the
proofs in~\cite{abhhs:noncommutative,a:repCalg} to which we refer for
details.

Let us start by studying the case when $\lambda\neq 0$, and let
$\xv=(e^{i\vphi},e^{i\vphit})$ be a vector such that $\vphi-\vphit=\pm
2\beta$. Then it is easy to calculate that $s(\xv)=(e^{i(\vphi\pm
  2\beta)},e^{i\vphi})$. Thus, the maps $s$ preserves the condition
(\ref{eq:philambda}) provided that we start with a vector fulf\/illing
the condition.

This implies that if $G_\L$ has a loop (i.e. a directed cycle) on $k$
vertices, then we must have that $s^k(\xv_i)=\xv_i$ for some $i\in
V$. Moreover, it is a trivial fact that every directed graph of a
unitary matrix must have a loop.  Hence, $\beta$ must be such that
$e^{i2k\beta}=1$ for some integer $k>0$.  Since the map $s$ is
invertible, given any $\xv_i$ in the graph uniquely determines $\xv_j$
for all other vertices in the graph. Hence, we can partition the
vertices into subsets $V_1,\ldots,V_k$ such that $\xv_i=\xv_j$ if and
only if $i,j\in V_l$ for some $l\in\{1,\ldots,k\}$. It follows that
all edges of $G_\L$ are of the form $(i,j)$ with $i\in V_l$ and $j\in
V_{l+1}$ (where we identify $k+1\equiv 1$). Thus, we can permute the
vertices to bring the matrix $\L$ to the following form
\begin{align*}
  \L=
  \begin{pmatrix}
    0               &  \L_1  & 0                & \cdots & 0 \\
    0               &  0                & \L_2  & \cdots & 0 \\
    \vdots          & \vdots            & \ddots           & \ddots & \vdots\\
    0               & 0                 & \cdots           & 0 &  \L_{k-1} \\
    \L_k & 0                 & \cdots           & 0 & 0
  \end{pmatrix}
\end{align*}
with $\L_i$ being unitary matrices for $i=1,\ldots,k$. Moreover, there
exists a unitary matrix such that $U^\dagger\L U$ is of the above form
but each $\L_i$ is diagonal. This means that the directed graph of
$U^\dagger\L U$ is a direct sum of $k$ loops (this also holds for $T$
since $T=\Dt\L^\dagger$), and each of these loops correspond to an
irreducible representation. However, to calculate the representation
graph, we must go to the basis in which $\L$ is diagonal. The matrix
corresponding to a single loop on $n$ vertices has the $n$ roots of
unity as eigenvalues. Therefore, in the basis in which $\L$ is
diagonal, the directed graph of $T$ will be the representation graph.
It is easy to see that the matrices $D$ and $\Dt$ will act as shift
operators on the eigenvectors of $\Lambda$, which implies that
$T=\Dt\L^\dagger$ will also act as a shift operator in this
basis. Thus, the directed graph of $T$ will be a single loop. We
conclude that this representation is precisely the Fuzzy torus
representation presented above.

\begin{figure}[h]
  \centering
  \includegraphics[width=7cm]{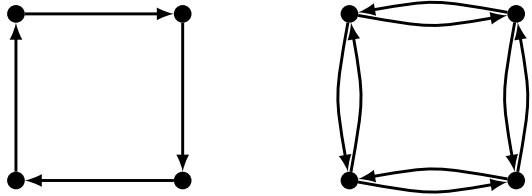}
  \caption{The two dif\/ferent types of representation graphs for
    irreducible unitary representations when $\mu_1=\cdots=\mu_4=2$.}
  \label{fig:unitary4d}
\end{figure}

Let us turn to the case when $\lambda=0$, i.e.~$\mu=2$. In this case,
there is no restriction on $\vphi-\vphit$, but instead one notes that
$s^4(\xv)=\xv$ for any $\xv\in\complex^2$. Thus, the vertices of
$G_\L$ can be split into four disjoint subsets, and we conclude that
all irreducible representations are 4-dimensional. However, since
$\vphi-\vphit$ does not have to be related to $\beta$, the action of
$T$ on the eigenbasis of $\Lambda$ will not simply be a
shift. Therefore, the representation graph will be one of the two in
Fig.~\ref{fig:unitary4d}.

When $\lambda=2$ ($\mu=0$) then $\vphi=\vphit$ and we see that any
vector of the form $(e^{i\vphi},e^{i\vphi})$ is a f\/ixpoint of
$s$. Hence, all irreducible representations are 1-dimensional. This
agrees with the result in~\cite{hs:repYangMillsAlgebras} which states
that, when $\Spec(\A)=\{0\}$, all irreducible f\/inite-dimensional
representations of~$\U_d(\A)$ are 1-dimensional.

\begin{proposition}
  Assume that $\mu\neq 2$ and let $\A$ be a \DMSA\ with
  $\Spec(\A)=\{\mu\}$.  If $\phi$ is an $n$-dimensional irreducible
  unitary representation of $\U_4(\A)$ then $\phi$ is equivalent to a
  representation~$\phi'$ with $\phi'(t_1)=e^{i\theta}g$ and
  $\phi'(t_2)=e^{i\theta'}h$ for some
  $\theta,\theta'\in\reals$. Moreover, there exists a $\beta\in\reals$
  such that $\mu=4\sin^2(\beta)$ and $e^{i2n\beta}=1$.
\end{proposition}

\begin{proposition}
  Let $\A$ be a \DMSA\ with $\Spec(\A)=\{2\}$ and let $\phi$ be an
  irreducible unitary representation of $\U_4(\A)$. Then $\phi$ is
  $4$-dimensional and the representation graph of $\phi$ is one of the
  two in Fig.~{\rm \ref{fig:unitary4d}}.
\end{proposition}

 The graph to the right in Fig.~\ref{fig:unitary4d} is the
Cartesian product of two representation graphs corresponding to
2-dimensional representations def\/ined by (\ref{eq:sphereNormalRep1})
and (\ref{eq:sphereNormalRep2}). However, one can check that there are
4-dimensional unitary representations that can not be written as a
tensor product of two such representations.

\subsection[Representations induced by $\sl_3$]{Representations induced by $\boldsymbol{\sl_3}$}

  We will now present a \DMSA\ $\A$ constructed from $\sl_3$,
whose representations give rise to normal representations
of $\U_4(\A)$. The vertices of the representation graph will be the
weight diagram of the $\sl_3$-representation.

\begin{figure}[h]
  \centering
  \includegraphics[width=5cm]{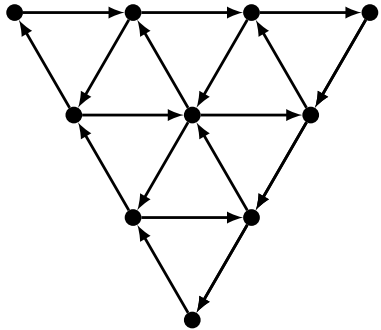}
  \caption{The representation graph corresponding to the $\{3,0\}$ representation
    of $\sl_3$.}
\label{fig:30rep}
\end{figure}

Let $\alpha$ and $\beta$ be the simple roots of $\sl_3$. By setting
\begin{gather*}
   t_1 = e^{i\theta}\big(h_\alpha+e^{i\pi/3}h_\beta\big),\qquad
  s_1 = e^{-i\theta}\big(h_\alpha+e^{-i\pi/3}h_\beta\big),\\
   t_2 = e^{i\theta'}\big(e_\alpha+e^{i\vphi_1}e_\beta+e^{i\vphi_2}e_{-\alpha-\beta}\big),\qquad
   s_2 = e^{-i\theta'}\big(e_{-\alpha}+e^{-i\vphi_1}e_{-\beta}+e^{-i\vphi_2}e_{\alpha+\beta}\big)
\end{gather*}
we obtain a \DMSA\ with $[t_1,s_1]=[t_2,s_2]=0$ and
\begin{gather*}
   \Ccom{t_1}{t_2}{s_2} = \frac{3}{2}l^2t_1,\qquad
  \Ccom{s_1}{t_2}{s_2} = \frac{3}{2}l^2s_1,\\
   \Ccom{t_2}{t_1}{s_1} = \frac{3}{4}l^4t_2,\qquad
  \Ccom{s_2}{t_1}{s_1} = \frac{3}{4}l^4s_2.
\end{gather*}
In the current convention, a compact real form of $\sl_3$ is provided
by
\begin{align*}
  ih_\alpha, \ ih_\beta, \ \eap, \ \eam, \ \ebp, \ \ebm, \
  e_{\alpha+\beta}^+, \ e_{\alpha+\beta}^-
\end{align*}
as def\/ined in Lemma \ref{lemma:ealphaplus}. Hence, any representation
is equivalent to one where these elements are represented by
\emph{anti-Hermitian} matrices, which implies that
$\phi(e_\gamma)^\dagger=\phi(e_{-\gamma})$ and
$\phi(h_\gamma)^\dagger=\phi(h_\gamma)$ for
$\gamma=\alpha,\beta,\alpha+\beta$. It follows that
$\phi(t_1)^\dagger=\phi(s_1)$ and $\phi(t_2)^\dagger=\phi(s_2)$.  The
weight diagram of a~representation is usually presented as vectors
with respect to an orthonormal basis of the Cartan subalgebra. In
$\sl_3$ we can construct an orthonormal basis by setting
\begin{gather*}
  h_1 = \frac{1}{l}h_\alpha,\qquad
  h_2 = \frac{1}{l\sqrt{3}}(h_\alpha+2h_\beta),
\end{gather*}
and from this we calculate that
\begin{gather*}
  t_1 = \frac{l\sqrt{3}}{2}e^{i(\theta-\pi/6)}(h_1+ih_2).
\end{gather*}
Hence, the eigenvalues of $\phi(t_1)$ will be the weights of the
(scaled and rotated) weight diagram of the representation $\phi$.  As
an example, let us study the representations of the kind $\{n,0\}$,
i.e.\ representations of highest weight $nw_1$, where $w_1$, $w_2$ are the
fundamental weights. These representations have dimension
$(n+1)(n+2)/2$ and all weights have multiplicity one. Therefore, in a~representation $\phi$, where the elements of the Cartan subalgebra are
diagonal, the representation graph is given by the directed graph of
$T=\phi(t_2)$. Since $t_2$ is a linear combination of $e_\alpha$,
$e_\beta$ and~$e_{-\alpha-\beta}$ we can construct the representation
graph by drawing arrows in the direction of these roots in the weight
diagram, see Fig.~\ref{fig:30rep}.

\subsection{Representations induced from a Clif\/ford algebra}

\begin{figure}[h]
  \centering
  \includegraphics[width=3cm]{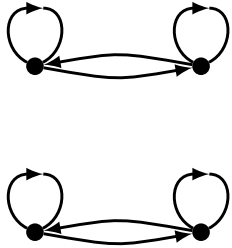}
  \caption{The representation graph of an irreducible representation
    related to the Clif\/ford algebra~$\Cl_{4,0}$.}
  \label{fig:cliffordrep}
\end{figure}

 As we have already noted, Clif\/ford algebras satisfy the
relations in the enveloping algebra and therefore, any
Clif\/ford algebra representation will give a representation of
$\U_d(\A)$. As an example, we consider the unique (up to equivalence)
irreducible representation of $\Cl_{4,0}$ given by
$e_i=\sigma_1\otimes\sigma_i$, for $i=1,2,3$, and
$e_4=\sigma_3\otimes\mid_2$, where
\begin{gather*}
  \sigma_1=\twomatrix{0}{1}{1}{0},\qquad
  \sigma_2=\twomatrix{0}{-i}{i}{0},\qquad
  \sigma_3=\twomatrix{1}{0}{0}{-1}.
\end{gather*}
The \emph{Hermitian} matrices $e_1$, $e_2$, $e_3$, $e_4$ satisfy the relations
$e_ie_j+e_je_i=2\delta_{ij}\mid_4$. In this example, no combination of
the form $e_k+ie_l$ will be diagonalizable which, in particular, means
that such a combination is never a normal matrix. However, as this is
an irreducible representation, any Jordan basis of $\Lambda=e_1+ie_2$,
is a Jordan basis of $\Lambda$ with respect to $T=e_3+ie_4$. The
Jordan normal form of $\Lambda$, as well as the matrix form of $T$ in
this basis, are easily computed to be
\begin{gather*}
  P^{-1}\Lambda P =
  \begin{pmatrix}
    0 & 1 & 0 & 0 \\
    0 & 0 & 0 & 0 \\
    0 & 0 & 0 & 1 \\
    0 & 0 & 0 & 0
  \end{pmatrix},
  \qquad
  P^{-1}T P =
  \begin{pmatrix}
    1+i & 0 & 1 & 0 \\
    0 & -1-i & 0 & -1 \\
    -2i & 0 & -1-i & 0 \\
    0 & 2i & 0 & 1+i
  \end{pmatrix}.
\end{gather*}
Hence, a representation graph $G_\phi=(G,\lambda)$ is given as in
Fig.~\ref{fig:cliffordrep} together with
\mbox{$\lambda:\{1,2,3,4\}\to\complex$} def\/ined by $\lambda(i)=0$ for all
$i$. We note that even though the graph is disconnected the
representation is irreducible.

\section{Summary}

 Motivated by several examples, in which double commutator
matrix equations arise, we have considered the relations
\begin{gather}
  \sum_{j=1}^d\Ccom{x_i}{x_j}{x_j}=\mu_i x_i\label{eq:xdoublecom}
\end{gather}
in a general (Lie) algebraic setting. Some examples can easily be
constructed from subsets of semi-simple Lie algebras.  Via the Diamond
lemma we can show that it is consistent to impose relations
(\ref{eq:xdoublecom}) in a free associative algebra, and a basis for
the corresponding enveloping algebra was computed.

In contrast to the case when $\mu_i=0$ for $i=1,\ldots,d$ (in which
case all irreducible f\/inite-dimensional representations are one
dimensional \cite{hs:repYangMillsAlgebras}), the representation theory
for arbitrary $\mu_i$'s has a rich structure. We have considered the
case when $d\leq 4$ in detail and introduced the \emph{representation
  graph}, which encodes the structure of a f\/inite-dimensional
representation in terms of a directed graph. The connectivity of the
graph provides information on the irreducible components of the
representation, and the tensor product can (generically) be described
by the Cartesian product of graphs.  All unitary representations when
$\Spec(\A)=\{\mu\}$ were then classif\/ied, and it was shown that
essentially all such representations are equivalent to the Fuzzy torus
algebra.  Several other examples were provided to demonstrate that the
representation theory is non-trivial. A particular feature,
that in general distinguishes the representation theory from that of
Lie algebras, is that the tensor product of two irreducible
representations can again be irreducible.

While we think relations (\ref{eq:xdoublecom}) are interesting in
themselves -- as a class of algebras containing Clif\/ford algebras and
semi-simple Lie algebras -- let us end by noting that we expect matrix
sequences corresponding to surfaces of genus $g\geq 2$ to exist, even
for $d\leq 4$ and $\mu_1=\cdots=\mu_4$.

\subsection*{Acknowledgement}

 {\sloppy We would like to thank the Marie Curie Research Training
Network ENIGMA and the Swedish Research Council, as well as the IHES,
the Sonderforschungsbereich ``Raum-Zeit-Materie'' (SFB647) and ETH
Z\"urich, for f\/inancial support respectively hospitality -- and Martin
Bordemann for many discussions and collaboration on related topics.

}

\pdfbookmark[1]{References}{ref}
\LastPageEnding


\begin{thebibliography}{99}

\footnotesize\itemsep=0pt

\bibitem{a:phdthesis}
 Arnlind J.,
 Graph techniques for matrix equations and eigenvalue dynamics,
 PhD Thesis, Royal Institute of Technology, 2008.

\bibitem{a:repCalg}
Arnlind J.,
Representation theory of $C$-algebras for a higher-order class of spheres and tori,
\href{http://dx.doi.org/10.1063/1.2913523}{{\em J. Math. Phys.}} {\bf 49} (2008), 053502, 13~pages,
\href{http://arxiv.org/abs/0711.2943}{arXiv:0711.2943}.

\bibitem{abhhs:noncommutative}
 Arnlind J.,  Bordemann M.,  Hofer L.,  Hoppe J.,   Shimada H.,
Noncommutative Riemann surfaces by embeddings in $\mathbb{R}^3$,
\href{http://dx.doi.org/10.1007/s00220-009-0766-8}{{\em Comm. Math. Phys.}} {\bf  288} (2009), 403--429,
\href{http://arxiv.org/abs/0711.2588}{arXiv:0711.2588}.

\bibitem{aht:spinning}
 Arnlind J.,  Hoppe J.,  Theisen S.,
 Spinning membranes,
\href{http://dx.doi.org/10.1016/j.physletb.2004.08.026}{{\em Phys. Lett. B}} {\bf 599} (2004), 118--128.


\bibitem{bdv:inhomyangmills}
 Berger R.,  Dubois-Violette M.,
 Inhomogeneous Yang--Mills algebras,
\href{http://dx.doi.org/10.1007/s11005-006-0075-5}{{\em Lett. Math. Phys.}} {\bf  76} (2006), 65--75,
\href{http://arxiv.org/abs/math.QA/0511521}{math.QA/0511521}.

\bibitem{b:diamondLemma}
 Bergman G.M.,
 The diamond lemma for ring theory,
\href{http://dx.doi.org/10.1016/0001-8708(78)90010-5}{{\em Adv. in Math.}} {\bf 29} (1978), 178--218.

\bibitem{b:threeDimLieAlgebras}
 Bianchi L.,
 Sugli spazi a tre dimensioni che ammettono un gruppo continuo di   movimenti,
 {\em Mem. Soc. Ital. delle Scienze (3)} {\bf  11} (1898), 267--352.

\bibitem{cdv:yangMillsAlgebra}
 Connes A.,  Dubois-Violette M.,
Yang--Mills algebra,
\href{http://dx.doi.org/10.1023/A:1020733628744}{{\em Lett. Math. Phys.}} {\bf 61} (2002), 149--158,
\href{http://arxiv.org/abs/math.QA/0206205}{math.QA/0206205}.

\bibitem{ffz:trigonometric}
 Fairlie D.B., Fletcher P., Zachos C.K.,
Trigonometric structure constants for new inf\/inite-dimensional algeb\-ras,
\href{http://dx.doi.org/10.1016/0370-2693(89)91418-4}{{\em Phys. Lett. B}} {\bf 218} (1989), 203--206.

\bibitem{ht:connectednessGraph}
Harary F., Trauth C.A. Jr.,
 Connectedness of products of two directed graphs.
\href{http://dx.doi.org/10.1137/0114024}{{\em SIAM J. Appl. Math.}} {\bf 14} (1966), 250--254.

\bibitem{hs:repYangMillsAlgebras}
 Herscovich E., Solotar A.,
Representation theory of Yang--Mills algebras,
\href{http://arxiv.org/abs/0807.3974}{arXiv:0807.3974}.

\bibitem{h:phdthesis}
 Hoppe J.,
 Quantum theory of a massless relativistic surface and a two-dimensional bound state problem,
PhD Thesis, MIT, 
 1982, available at
\url{http://dspace.mit.edu/handle/1721.1/15717}.

\bibitem{h:diffeomorphism}
 Hoppe J.,
Dif\/feomorphism groups, quantization, and ${\rm SU}(\infty)$,
\href{http://dx.doi.org/10.1142/S0217751X89002235}{{\em Internat. J. Modern Phys. A}} {\bf 4} (1989), 5235--5248.

\bibitem{h:someClassicalSolutions}
 Hoppe J.,
Some classical solutions of membrane matrix model equations,
\href{http://arxiv.org/abs/hep-th/9702169}{hep-th/9702169}.

\bibitem{l:completeMinimalSurfaces}
Lawson H.B. Jr.,
Complete minimal surfaces in $S\sp{3}$,
\href{http://dx.doi.org/10.2307/1970625}{{\em Ann. of Math. (2)}} {\bf 92} (1970), 335--374.

\bibitem{n:lecturenoncgauge}
Nekrasov N.,
Lectures on open strings, and noncommutative gauge theories,
\href{http://arxiv.org/abs/hep-th/0203109}{hep-th/0203109}.

\bibitem{s:graphMultiplication}
Sabidussi G.,
Graph multiplication,
\href{http://dx.doi.org/10.1007/BF01162967}{{\em Math. Z.}} {\bf  72} (1959/1960), 446--457.

\end{thebibliography}
\end{document}